\documentclass[reqno,11pt]{amsart}
\usepackage{latexsym}
\usepackage{amsfonts}
\usepackage{amssymb}
\usepackage{mathrsfs}
\usepackage[all]{xy}
\usepackage{bbm}
\usepackage{color}
\usepackage[mathscr]{euscript}

\usepackage{hyperref}

\newtheorem{theorem}{Theorem}[section]
\newtheorem{corollary}[theorem]{Corollary}
\newtheorem{lemma}[theorem]{Lemma}

\newtheorem{proposition}[theorem]{Proposition}
\theoremstyle{definition}
\newtheorem{definition}[theorem]{Definition}

\numberwithin{equation}{section}
\newtheorem{question}[theorem]{Question}

\newcommand{\TT}{\mathbb{T}}
\newcommand{\RR}{\mathbb{R}}
\newcommand{\NN}{\mathbb{N}}

\newcommand{\ZZ}{\mathbb{Z}}
\newcommand{\be}{\mathbf{e}}
\newcommand{\bff}{\mathbf{f}}
\newcommand{\cA}{\mathcal{A}}

\newcommand{\cD}{\mathcal{D}}

\newcommand{\cW}{\mathcal{W}}

\newcommand{\Bad}{\mathbf{Bad}}
\newcommand{\bx}{\mathbf{x}}

\newcommand{\bv}{\mathbf{v}}
\newcommand{\br}{\mathbf{r}}
\newcommand{\bs}{\mathbf{s}}
\newcommand{\bw}{\mathbf{w}}
\newcommand{\bp}{\mathbf{p}}
\newcommand{\bq}{\mathbf{q}}

\newcommand{\diag}{\mathrm{diag}}
\newcommand{\Span}{\mathrm{Span}}
\newcommand{\dist}{\mathrm{dist}}
\newcommand{\supp}{\mathrm{supp}}

\newcommand{\SL}{\mathrm{SL}}

\newcommand{\id}{\mathrm{id}}

\newcommand{\la}{\left\langle}
\newcommand{\ra}{\right\rangle}

\newcommand {\ignore}[1]  {}

\newif\ifdraft\drafttrue

\draftfalse

\renewcommand{\emptyset}{\varnothing}
\renewcommand{\setminus}{\smallsetminus}

\begin{document}
	
	\title{$\Bad(\br;\bs)$ is Hyperplane Absolute Winning}
\date{}
\author{Chengyang Wu}
\address{Department of Mathematics, University of Chicago, Chicago, Illinois, 60637, U.S.
}
\email{chengyangwu1999@gmail.com}

	\begin{abstract}
    Given an $m$-dimensional weight $\br$ and an $n$-dimensional weight $\bs$, we prove that the set of $(\br;\bs)$-badly approximable $m\times n$ matrices is hyperplane absolute winning on $\RR^{m\times n}$. This fully answers a question \cite[Question 8.2 (iii)]{Kl} of D. Kleinbock in 1998.
	\end{abstract}

    \maketitle

    \section{Introduction}

    \subsection{Weighted badly approximable matrices}
     For $a\in\NN$, let $\cW_a$ denote the set of all $a$-dimensional weights: \begin{equation}
        \cW_a:=\left\{
        \bw=(w_1,\cdots,w_a)\in\RR^a:w_1,\cdots,w_a>0,\;\sum_{i=1}^aw_i=1
        \right\}.
    \end{equation}
    Here, for the sake of simplicity, we require that all weight components are positive (see Section \ref{SS:zeroweight} for the extension). Following the notation in \cite{Kl}, for $\bw\in\cW_a$, we define the $\bw$-\textsl{quasi-norm} $\|\cdot\|_{\bw}$ on $\RR^a$ by \begin{equation}
        \|\bx\|_{\bw}:=\max_{1\leq i\leq a}|x_i|^{1/w_i}.
    \end{equation}

    Weighted Diophantine approximation shares many important features with its classical unweighted counterpart. For example, the fundamental Dirichlet's theorem easily extends to the weighted case: for $(\br;\bs)\in\cW_m\times \cW_n$, any $X\in\RR^{m\times n}$, and any $Q>1$, there exist $\bp\in\mathbb{Z}^m$ and $\bq\in\mathbb{Z}^n\setminus\{0\}$ such that
$$
    |(X\bq+\bp)_i|\le Q^{-r_i}\quad (1\le i\le m),
    \qquad |\bq_j|\le Q^{s_j}\quad (1\le j\le n).
$$

    The focus of this paper is the set of weighted badly approximable matrices. Let $m,n$ be two positive integers. For $(\br;\bs)\in\cW_m\times \cW_n$, an $m\times n$ real matrix $X$ is $(\br;\bs)$-\textsl{badly approximable}, if \begin{equation}\label{E:Bad}
        \inf_{\substack{\bp\in\ZZ^m\\ \bq\in\ZZ^n\setminus\{0\}}}
        \|X\bq+\bp\|_{\br}\cdot \|\bq\|_{\bs}>0.
    \end{equation}
    The set of $(\br;\bs)$-badly approximable matrices is denoted by $\Bad(\br;\bs)$. It is of particular interest to investigate its qualitative and quantitative properties.

    In view of Khintchine-Groshev's theorem \cite{Gro38}, it is not difficult to see that $\Bad(\br;\bs)$ is null with respect to the standard Lebesgue measure on $\RR^{m\times n}$. On the other hand, it was proved by D. Kleinbock and B. Weiss in \cite{KW08,KW10,KW13} that $\Bad(\br;\bs)$ is \textsl{thick} in $\RR^{m\times n}$, that is, the set has full Hausdorff dimension in any nonempty open subset in the ambient space.

    However, to the best of our knowledge, for two general pairs $(\br_1;\bs_1),(\br_2;\bs_2)\in\cW_m\times \cW_n$, it was not even known whether \begin{equation}\label{E:inter}
        \Bad(\br_1;\bs_1)\cap \Bad(\br_2;\bs_2)\neq \varnothing.
    \end{equation}
    This can be viewed as a matrix version of Schmidt's conjecture \cite{Sch83}, which concerns the case $(m,n)=(2,1)$. Over the past fifteen years, Schmidt's conjecture was proved in \cite{BPV}, and its higher-dimensional analogue for vectors was established in \cite{Be}.

    An alternative framework to show the nonempty intersection \eqref{E:inter} is provided by game theory. This was first developed by Schmidt in his seminal work \cite{Sch66}, where a game and its winning sets are introduced. Here we list two basic properties of winning sets, and readers may refer to Section \ref{SS:Schmidt} for details: \begin{itemize}
        \item winning subsets of Euclidean spaces are thick;
        \item a countable intersection of winning sets is again winning.
    \end{itemize}

    Ever since, Schmidt's game and its variants have been powerful tools in proving thickness of a finite or countable intersection. For example, one may easily deduce \eqref{E:inter} from the winning property of $\Bad(\br;\bs)$ for two pairs of weights $(\br;\bs)$. In \cite{Kl}, Kleinbock proposed the following general question: 
    
    \begin{question}[{\cite[Question 8.2 (iii)]{Kl}}]\label{Q:Kl}
        For arbitrary choice of $\br$ and $\bs$, is $\Bad(\br;\bs)$ a winning set of $\RR^{m\times n}$?
    \end{question}
    
    Regarding Question \ref{Q:Kl}, the unweighted case was proved in \cite{BFKRW12,BFS13,Sch66,Sch69}, and the weighted cases for vectors were established in \cite{An16,BNY21,GY,NS}. So far, no weighted case for general matrices was resolved.

    This paper answers Question \ref{Q:Kl} in full generality.
    
    \begin{theorem}\label{T:main}
        For any $(\br;\bs)\in\cW_m\times\cW_n$, the set $\Bad(\br;\bs)$ is HAW on $\RR^{m\times n}$.
    \end{theorem}

    Here HAW stands for \textsl{hyperplane absolute winning} (see Section \ref{SS:Schmidt} for definitions and properties). In particular, the following corollary is immediate: 

    \begin{corollary}\label{C:main}
        Let $(\br_k;\bs_k)_{k\in\NN}$ be a sequence of pairs of weights in $\cW_m\times\cW_n$. Then the intersection $$
        \bigcap_{k\in\NN}\Bad(\br_k;\bs_k)
        $$
        is thick in $\RR^{m\times n}$.
    \end{corollary}

    \subsection{Relationship to homogeneous dynamics} As is well-known, badly approximable matrices correspond to certain bounded orbits on homogeneous spaces. Throughout this paper, set $d=m+n$. Fix $b>1$. For $(\br;\bs)\in\cW_m\times \cW_n$, we consider a one-parameter diagonal subgroup $F^{(\br;\bs)}=\{g_t^{(\br;\,\bs)}:t\in\RR\}$ in $G=\SL_d(\RR)$ as \begin{equation}
        g_t^{(\br;\,\bs)}=\diag\left((b^{r_it})_{1\leq i\leq m},\;(b^{-s_jt})_{1\leq j\leq n}\right).
    \end{equation}
    We also consider a unipotent subgroup $U_0=\{u_X:X\in \RR^{m\times n}\}$ in $G$, where \begin{equation}
        u_X=\begin{pmatrix}
            I_m&X\\ 0& I_n
        \end{pmatrix}.
    \end{equation}
    Then Dani-Kleinbock correspondence (see \cite{Dani85,Kl}) states that $X$ is $(\br;\bs)$-badly approximable if and only if the trajectory $\{g_t^{(\br;\,\bs)}u_X\Gamma:t\geq 0\}$ is bounded in $G/\Gamma$ where $\Gamma=\SL_d(\ZZ)$.

    Consider the diagonal subsemigroup in $G$:$$
    A^+=\left\{
    \diag((b^{u_i})_{1\leq i\leq m},(b^{-v_j})_{1\leq j\leq n}):u_i,v_j\geq 0,\;\sum_{i=1}^mu_i=\sum_{j=1}^nv_j
    \right\}.
    $$
    A one-parameter subsemigroup in $A^+$ is called a \textsl{diagonal ray}. For example, the semigroup \begin{equation}\label{E:diagonal-ray}
    (F^{(\br;\bs)})^+:=\{g_t^{(\br;\,\bs)}:t\geq 0\}
    \end{equation}
    is a diagonal ray; conversely, any diagonal ray without eigenvalue $1$ is of the form \eqref{E:diagonal-ray}.

    In the language of homogeneous dynamics, we reformulate Theorem \ref{T:main} and Corollary \ref{C:main} as corresponding dynamical statements on the $mn$-dimensional subtorus in $G/\Gamma$: \begin{equation}
        \TT_0:=\{u_X\Gamma:X\in\RR^{m\times n}\}.
    \end{equation}  

    \begin{theorem}\label{T:dyn}
        For any diagonal ray $F^+$ in $A^+$, the set $$
        B(F^+):=\{p\in\TT_0:F^+p\text{ is bounded}\}
        $$
        is HAW on $\TT_0$.
    \end{theorem}

    \begin{corollary}
        Let $(F^+_k)_{k\in\NN}$ be a sequence of diagonal rays in $A^+$. Then the intersection $$
        \bigcap_{k\in\NN}B(F_k^+)=\{p\in\TT_0:F_k^+p\text{ is bounded for all }k\in\NN\}
        $$
        is thick in $\TT_0$.
    \end{corollary}

    \subsection{Organization of the paper and strategy of the proof} The whole paper is organized as follows. In Section \ref{S:pre}, we recall some preliminary knowledge on Schmidt's games and two of its variants: the hyperplane absolute game and Cantor potential game, especially their winning criterion (see Theorem \ref{T:BNY21}). We also recall the quantitative non-divergence result (see Proposition \ref{P:BNY21}) for unipotent actions with respect to an Ahlfors regular, absolutely decaying measure. These useful tools have been employed in the previous papers \cite{BNY21,BNY22}.
    
    Sections \ref{S:global-estimate} and \ref{S:local-estimate} consist of the main body of our paper. The essential technicality is to introduce the family of diagonal matrices $\{D_t^{(k,\ell)}:t\in\RR\}$ for $(k,\ell)\in\NN_0^2$, and to estimate the supremum norm of the translated $r$-vector $\bigwedge^r(D^{(k,\ell)}_tu_X)\xi$ over a ball. These will exclude the possibility (ii) in Proposition \ref{P:BNY21} and hence yield the desired measure control of dangerous sets as (i) in Proposition \ref{P:BNY21}.

    The estimate for the vector norm of $\bigwedge^r(D^{(k,\ell)}_tu_X)\xi$ is divided into the global and local parts. In Section \ref{S:global-estimate}, the global estimate is deduced from the polynomial nature of the map \begin{equation}\label{E:poly}
        X\mapsto \bigwedge^r(D^{(k,\ell)}_tu_X)\xi.
    \end{equation}
    In particular, when the ambient ball $B$ is not too small, the supremum norm of \eqref{E:poly} over $B$ at certain time $t$ has a uniform positive lower bound.

    The case is more involved when the ambient ball $B$ is comparatively small. In Section \ref{S:local-estimate}, the local estimate relies on the key observation: Lemma \ref{L:dicho} in exterior algebra. Basically, suppose that the supremum norm of \eqref{E:poly} at certain time $t$ is small over a local ball $2B$.
    For $X_0\in B$, by considering the action of $D^{(k,\ell)}_tu_{X-X_0}(D^{(k,\ell)}_t)^{-1}$ (where $\|X-X_0\|\leq\rho(B)$) on the $r$-vector \begin{equation}
    \zeta:=\bigwedge^r(D^{(k,\ell)}_tu_{X_0})\xi,
    \end{equation} 
    we see that a spatial segment becomes a long shear. In particular, when this shear happens along a rank-$1$ direction, it produces two different diagonal flows, under the action of at least one of which the vector norm of $\zeta$ gets shrunk. Then this dangerous point $X_0$ should go to the historical account at some previous game stage.

    Therefore, we conclude that as long as the current ball $B$ is not covered by the historical dangerous sets, the supremum norm of \eqref{E:poly} at certain time $t$ has a uniform positive lower bound over a local ball $2B$. Then applying Proposition \ref{P:BNY21} gives the local measure control \eqref{E:measure-control} over dangerous sets of certain scales. 

    In Section \ref{S:Proofs}, we verify that this local measure control over dangerous sets guarantees a legal collection of balls for Alice to delete in the Cantor potential game. Hence Alice has a winning strategy with the target set $\Bad(\br;\bs)$, which proves Theorem \ref{T:main}. The extension of Theorem \ref{T:main} to weights with zero components follows from a simple property of winning sets. Finally, we prove that Theorem \ref{T:dyn} is equivalent to this extended result.

    \section{Preliminaries}\label{S:pre}

    \subsection{Schmidt's game and its variants}\label{SS:Schmidt} In this subsection, we recall Schmidt's game and two of its variants: the hyperplane absolute game and Cantor potential game.

    \begin{definition}
Schmidt's game \cite{Sch66} involves two parameters $\alpha,\beta\in(0,1)$  and is played by two players Alice and Bob on a complete metric space $(X,\dist)$ with a target set $S\subset X$. Bob starts the game by choosing a closed ball $B_0$ of radius $\rho_0$. After Bob chooses a closed ball $B_q$ of radius $\rho_q$, Alice chooses a ball $A_q\subseteq B_q$ of radius ${\rho}'_q=\alpha {\rho}_q$, and then Bob chooses a ball $B_{q+1}\subseteq A_q$ of radius ${\rho}_{q+1}=\beta {\rho}'_{q}$. This implies that the balls are nested:
$$
	B_0 \supset A_0 \supset B_1 \supset\cdots.$$
Alice wins the game if the unique point in $\bigcap_{q=0}^\infty A_q=\bigcap_{q=0}^\infty B_q$ belongs to $S$. The set $S$ is \textsl{$(\alpha,\beta)$-winning} if Alice has a winning strategy, is \textsl{$\alpha$-winning} if it is $(\alpha,\beta)$-winning for any $\beta\in(0,1)$, and is \textsl{winning} if it is $\alpha$-winning for some $\alpha\in (0,1)$.
\end{definition}

Regarding Schmidt's game and its winning sets, Schmidt \cite{Sch66} proved that:	
\begin{itemize}
	\item[(W1)]  winning subsets of Riemannian manifolds are thick;
    \item[(W2)]  a countable intersection of $\alpha$-winning sets is again $\alpha$-winning;
    \item[(W3)]  the image of a winning set under a bi-Lipschitz map is winning.
\end{itemize}
Schmidt's game has been a powerful tool for proving thickness of intersections of certain countable families of sets, see e.g.\ \cite{An13, An16,  BBFKW10, BFK11, BFKRW12}. Over the following decades, several useful and important variants of this game have been developed.


\begin{definition}
    Hyperplane absolute game \cite{BFKRW12} involves one parameter $\beta\in (0,\frac{1}{3})$ and is played by two players Alice and Bob on the Euclidean space $\RR^a$ with a target set $S \subset \RR^a$. Bob starts the game by choosing a closed ball $B_0$ of radius ${\rho}_0$.
For an affine hyperplane $L\subset\RR^a$ and $r>0$, we denote the $r$-neighborhood of $L$ by
$$L^{(r)}:=\{z\in\RR^a:\mathrm{dist}(z,L)\leq r\}.$$
After Bob chooses a closed ball $B_q$ of radius ${\rho}_q$, Alice chooses a hyperplane neighborhood $L_q^{(r_q)}$ with $r_q\le\beta {\rho}_q$, and then Bob chooses a closed ball $B_{q+1}\subset B_q\setminus L_q^{(r_q)}$ of radius ${\rho}_{q+1}\ge\beta {\rho}_q$. 

We say that Alice wins by default if at any time Bob has no legal move, or $\rho_q \not\to 0$ as $q\to +\infty$. Otherwise, Alice wins the game if and only if
$$\bigcap_{q=0}^{+\infty} B_q\cap S\ne\emptyset.$$ The set $S$ is \textsl{$\beta$-hyperplane absolute winning} (\textsl{$\beta$-HAW}) if Alice has a winning strategy, and is \textsl{hyperplane absolute winning} (\textsl{HAW}) if it is $\beta$-HAW for any $\beta\in(0,\frac{1}{3})$.
\end{definition}

\begin{lemma}[\cite{BFKRW12}]\label{L:HAW-Ra}
	\begin{itemize}
		\item[(i)] HAW subsets   are winning, and hence thick.
		\item[(ii)] A countable intersection of HAW subsets   is again HAW.
		\item[(iii)]   The image of an HAW set under a $C^1$ diffeomorphism $\RR^a \to \RR^a$ is HAW.
	\end{itemize}
\end{lemma}

\begin{definition}
    Cantor potential game \cite{BH17,BHNS18} involves two parameters $\alpha\geq 0$, $\beta\in(0,1)$,  and is played by two players Alice and Bob on a complete metric space $(X,\dist)$ with a target set $S\subset X$. Bob starts the game by choosing a closed ball $B_0$ of radius $\rho_0$ in $X$. After Bob chooses a closed ball $B_q$ of radius $\rho_0\beta^q$, Alice chooses, for each $i\geq 0$, a collection $\cA_{q+1,i}$ of at most $\beta^{-\alpha(i+1)}$ balls of radius $\rho_0\beta^{q+1+i}$. Then Bob chooses a ball $B_{q+1}$ of radius $\rho_0\beta^{q+1}$ such that $$
    B_{q+1}\subseteq B_q\setminus \bigcup_{0\leq q'\leq q}\bigcup_{A\in \cA_{q'+1,q-q'}}A.
    $$

    We say that Alice wins by default if at any time Bob has no legal move. Otherwise, Alice wins the game if and only if 
    $$\bigcap_{q=0}^{+\infty} B_q\cap S\ne\emptyset.$$ The set $S$ is \textsl{$(\alpha,\beta)$-Cantor winning} if Alice has a winning strategy, is \textsl{$\alpha$-Cantor winning} if it is $(\alpha,\beta)$-Cantor winning for any $\beta\in(0,1)$.
\end{definition}

It was shown in \cite{BH17} that on certain nice spaces (e.g. doubling spaces) the properties (W1)--(W3) still hold for Cantor winning sets. In this paper, we shall use another interesting property of Cantor winning sets. To do this, we need to recall the concept of \textsl{Ahlfors regular absolutely decaying measure}.

\begin{definition}\label{D:Ahlfors}
    Let $X$ be a metric space. Given $\alpha>0$, a Borel measure $\nu$ on $X$ is \textsl{$\alpha$-Ahlfors regular}, if there exist $A_0,r_0>0$ such that for any $x\in\supp(\nu)$ and any $r\in (0,r_0)$, 
$$A_0^{-1}r^{\alpha}\leq \nu(B(x,r))\leq A_0r^{\alpha}.$$
We say that $\nu$ is \textsl{Ahlfors regular} if it is $\alpha$-Ahlfors regular for some $\alpha>0$.
\end{definition}

\begin{definition}
    A Borel measure $\nu$ on $\RR^a$ is \textsl{absolutely decaying} if there exist $E_0,\delta,\rho_0>0$ such that for any $x\in\supp(\nu)$, any $\rho\in (0,\rho_0)$, every affine hyperplane $L\subseteq \RR^a$ and any $r>0$, one has $$
    \nu(B(x,\rho)\cap L^{(r)})\leq E_0(r/\rho)^{\delta}\cdot \nu(B(x,\rho)).
    $$
\end{definition}

We list here several external results for our later use:

\begin{lemma}[{\cite[Lemma 22]{BNY21}}]\label{L:BNY-cover}
    Let $\nu$ be an $\alpha$-Ahlfors regular measure on a metric space $X$, let $A_0,r_0>0$ be as in Definition \ref{D:Ahlfors}, and let $S\subseteq X$ be any subset. Suppose that $r\in (0,r_0)$ and $$\nu(S^{(r)}):=\nu\{x\in X:\dist(x,S)\leq r\}<+\infty.$$ Then there exists a cover of $S\cap\supp(\nu)$ by balls of radius $2r$ with cardinality at most $A_0\nu(S^{(r)})\cdot r^{-\alpha}$.
\end{lemma}

\begin{theorem}[{\cite[Proposition 11, Theorem 21]{BNY21}}]\label{T:BNY21}
\,\begin{enumerate}
    \item A Borel subset of $\RR^a$ that intersects the support of every nonzero compactly supported Ahlfors regular absolutely decaying measure is HAW.
    \item Let $X$ be the support of a nonzero $\alpha$-Ahlfors regular measure, and let $S\subseteq X$ be $\alpha'$-Cantor winning for some $0<\alpha'<\alpha$. Then $S\neq\varnothing$.
\end{enumerate}
\end{theorem}

    \subsection{Quantitative non-divergence arguments}
    The theory of quantitative non-diver-gence originates in Margulis's work \cite{Mar71} on the non-divergence of unipotent trajectories in the space of lattices, and Dani's subsequent recurrence estimates \cite{Dan79,Dan86}. In \cite{KM98}, Kleinbock and Margulis developed the $(C,\alpha)$-good function method and established the quantitative non-divergence estimates with power-law decay. Later, this framework was extended to certain non-Lebesgue measures including friendly measures \cite{KLW04}. These developments provide the quantitative non-divergence arguments employed in \cite{BNY21,BNY22}.

    To state the quantitative non-divergence result in the matrix setting considered here, we first introduce some notations of norms on exterior powers of Euclidean spaces and matrix spaces. Let $a\in\NN$ and let $\la\cdot,\cdot\ra$ denote the standard inner product on $\RR^a$. It induces a norm $\|\cdot\|$ on $\RR^a$ and hence a norm on the exterior powers $\bigwedge^r\RR^a\,(0\leq r\leq a)$. For a linear map $A$ on $\RR^a$, we write $\bigwedge^rA$ to be its induced map on $\bigwedge^r\RR^a$, and write the operator norm as $$
    \left\|\bigwedge^rA\right\|:=\sup_{\substack{\xi\in\bigwedge^r(\RR^a)\\\|\xi\|=1}}\left\|\bigwedge^rA(\xi)\right\|.
    $$
    We say that a nonzero $r$-vector $\xi$ is \textsl{decomposable} if it is a single wedge of $r$ vectors in $\RR^a$. For a lattice $\Lambda$ in $\RR^a$, write $$
    \lambda_1(\Lambda):=\min_{\bv\in\Lambda\setminus\{0\}}\|\bv\|.
    $$ 

    \begin{proposition}[{\cite[Proposition 26]{BNY21}}]\label{P:BNY21}
    Let $\nu$ be a nonzero, compactly supported, Ahlfors-regular, absolutely decaying measure on $\RR^{mn}$. Then there exist constants $C_0,\alpha_0,r_0>0$ such that for any ball $B$ centered in $\supp(\nu)$ with radius $\leq r_0$, any diagonal matrix $D\in\SL_d(\RR)$, and any $0<\rho\leq 1$, one of the following holds: \begin{enumerate}
            \item[(i)] for any $\epsilon>0$, \begin{equation}
                \nu\left\{
                X\in B:\lambda_1(Du_X\ZZ^d)<\epsilon
                \right\}\leq C_0(\epsilon/\rho)^{\alpha_0}\cdot \nu(3B).
            \end{equation}
            \item[(ii)] there exists some $1\leq r<d$ and a nonzero decomposable $r$-vector $\xi\in \bigwedge^r\ZZ^d$ with \begin{equation}
                \sup_{X\in B}\left\|
                \bigwedge^r(Du_X)\xi
                \right\|<\rho.
            \end{equation}
        \end{enumerate}    
    \end{proposition}

    \section{The Global Estimate}\label{S:global-estimate}

    In this section, we prove the global estimate for the supremum of $r$-vector norms, which serves as early steps in excluding the possibility (ii) in Proposition \ref{P:BNY21}.

    We first fix some notations. Let $(\be_1,\cdots,\be_m;\bff_1,\cdots,\bff_n)$ denote the standard ordered basis of $\RR^d=\RR^m\times \RR^n$. For a subset $I\subseteq \{1,\cdots,m\}$, write $$
    \be_I:=\begin{cases}
        \be_{i_1}\wedge\cdots\wedge \be_{i_a}, &\text{if }I=\{i_1<\cdots<i_a\} \text{ for some }a\geq 1;\\
        1, &\text{if }I=\varnothing.
    \end{cases}
    $$
    For a subset $J\subseteq\{1,\ldots,n\}$, we write $\bff_J$ similarly.
    
    Let $b>1$. For $(\br;\bs)\in\cW_m\times\cW_n$, we choose and fix some $i_*\in\{1,\cdots,m\}$ and $j_*\in\{1,\cdots,n\}$ such that  $$
    r_{i_*}=\max_{1\leq i\leq m}r_i>0,\qquad s_{j_*}=\max_{1\leq j\leq n}s_j>0.
    $$
    Write $$
    \lambda:=r_{i_*}+s_{j_*},\qquad 
    \mu:=\min\{r_1,\cdots,r_m,s_1,\cdots,s_n\}>0.
    $$
    For $k,\ell\in\NN_0$, we introduce a family of diagonal matrices $\{D_t^{(k,\ell)}:t\in\RR\}$ in $\SL_d(\RR)$ as \begin{equation}\label{E:Dkl-flow}
        D^{(k,\ell)}_t:=\diag\left(
        \left(b^{r_it-k1_{i=i_*}+\frac{k-\ell}{d}}\right)_{1\leq i\leq m},\;
        \left(b^{-s_jt+\ell1_{j=j_*}+\frac{k-\ell}{d}}\right)_{1\leq j\leq n}
        \right).
    \end{equation}
    The following proposition is the global estimate for the $r$-vector norms.

    \begin{proposition}\label{P:global-control}
        There exists an absolute constant $c_0>0$ (only depending on $m$ and $n$) with the following property. For any ball $B=B(X_0;\rho)$ in $\RR^{m\times n}$ with $0<\rho\leq 1$, any $k,\ell\in\NN_0$, any $t\geq 0$, and $1\leq r<d$, and any $\xi\in\bigwedge^r\ZZ^d\setminus\{0\}$, one has \begin{equation}\label{E:global-control}
            \sup_{X\in B}\left\|\bigwedge^r(D^{(k,\ell)}_tu_X)\xi\right\|\geq c_0b^{\mu t-(k+\ell)}\cdot \rho^{\min\{m,n\}}.
        \end{equation}
    \end{proposition}
    \begin{proof} Let $B=B(X_0;\rho)$, $k,\ell,t,r$, and $\xi$ be given, and write \begin{equation}
        P(Y):=\bigwedge^r(D^{(k,\ell)}_tu_{X_0+\rho Y})\xi.
    \end{equation}
    With respect to the fixed basis $\{\be_I\wedge\bff_J:|I|+|J|=r\}$ of $\bigwedge^r\RR^d$, write $P^{(I,J)}(Y)$ to be the $\be_I\wedge \bff_J$-coordinate of $P(Y)$. Note that each $P^{(I,J)}(Y)$ is a real-valued polynomial with degree $\leq d$ in $mn$ variables $Y=(Y_{ij})$. 
    
    Consider the finite-dimensional linear space consisting of all real-valued polynomial functions with degree $\leq d$ in $mn$ variables. There are two equivalent norms on this space: \begin{enumerate}
        \item[($\|\cdot\|_1$)] the supremum norm of function values over $B(0;1)\subseteq\RR^{m\times n}$;
        \item[($\|\cdot\|_2$)] the maximum norm of polynomial coefficients.
    \end{enumerate}
    It follows that for any pair of index sets $(I,J)$ with $|I|+|J|=r$, one has $$
    \text{LHS of \eqref{E:global-control}}\geq \|P^{(I,J)}(Y)\|_1\asymp \|P^{(I,J)}(Y)\|_2,
    $$
    where the implicit constant only depends on $m$ and $n$. To prove \eqref{E:global-control}, it suffices to show that: there exists a pair of index sets $(I_0,J_0)$ with $|I_0|+|J_0|=r$ such that \begin{equation}\label{E:I0J0}
        \|P^{(I_0,J_0)}(Y)\|_2\geq b^{\mu t-(k+\ell)}\cdot \rho^{\min\{m,n\}}.
    \end{equation}

    We first explicitly choose the pair $(I_0,J_0)$. Indeed, in the expression $$
    \xi=\sum_{|I|+|J|=r}\xi^{(I,J)}\be_I\wedge\bff_J \quad\text{where}\quad \xi^{(I,J)}\in\ZZ,
    $$
    we choose a pair $(I_1,J_1)$ with $|I_1|+|J_1|=r$ such that \begin{equation}\label{E:defJ1}
    \xi^{(I_1,J_1)}\neq 0,\quad\text{and}\quad |J_1|=\max\{|J|:\xi^{(I,J)}\neq 0\}.
    \end{equation}
    Set $h:=\min\{m-|I_1|,|J_1|\}$, and choose any $I_2\subseteq \{1,\cdots,m\}\setminus I_1$ and any $J_2\subseteq J_1$ with $|I_2|=|J_2|=h$. Then we put $I_0:=I_1\cup I_2$ and $J_0:=J_1\setminus J_2$. It is clear that $|I_0|+|J_0|=r$.

    We next verify that the chosen pair $(I_0,J_0)$ satisfies \eqref{E:I0J0}. Indeed, for the expression $$
    (\bigwedge^ru_{X_0})\xi=:\zeta=\sum_{|I'|+|J'|=r}\zeta^{(I',J')}\be_{I'}\wedge\bff_{J'} \quad\text{where}\quad \zeta^{(I',J')}\in\RR,
    $$
    by \eqref{E:defJ1}, we see that $\xi^{(I_1,J_1)}=\zeta^{(I_1,J_1)}$. Moreover, for the expression $$
    (\bigwedge^ru_{\rho Y})\zeta:=\theta=\sum_{|I|+|J|=r}\theta^{(I,J)}\be_I\wedge\bff_J \quad\text{where}\quad \theta^{(I,J)}\in\RR,$$ 
    we calculate that \begin{equation}\label{E:deftheta}
    \sum_{\substack{I'\subseteq I_0,\;J'\supseteq J_0\\ |I'|+|J'|=r}}(\pm 1)\zeta^{(I',J')}\cdot \rho^{|J'|-|J_0|}\det{Y_{I_0\setminus I',J'\setminus J_0}}=\theta^{(I_0,J_0)}.
    \end{equation}
    All monomials in $(Y_{ij})$ on the left-hand side of \eqref{E:deftheta} are different. In particular, it contains $$
    (\pm1 )\zeta^{(I_1,J_1)}\cdot \rho^{|J_1|-|J_0|}\det{Y_{I_0\setminus I_1,J_1\setminus J_0}}
    =(\pm 1)\xi^{(I_1,J_1)}\cdot \rho^h\det{Y_{I_2,J_2}}.
    $$
    Finally, the action of $\bigwedge^rD_t^{(k,\ell)}$ on $\be_{I_0}\wedge \bff_{J_0}$ is given by the scalar $$
    b^{(\sum_{i\in I_0}r_i-\sum_{j\in J_0}s_j)\cdot t-k1_{i_*\in I_0}+\ell 1_{j_*\in J_0}+(k-\ell)r/d}
    \geq b^{\mu t-(k+\ell)}.
    $$
    Therefore, we conclude that $$
    \|P^{(I_0,J_0)}(Y)\|_2\geq b^{\mu t-(k+\ell)}\|\theta^{(I_0,J_0)}\|_2\geq b^{\mu t-(k+\ell)}|\xi^{(I_1,J_1)}|\cdot \rho^h\geq b^{\mu t-(k+\ell)}\cdot \rho^{\min\{m,n\}}.
    $$
    This is as desired.
    \end{proof}

	\section{The Local Estimate}\label{S:local-estimate}
    In this section, we prove the local estimate for the supremum of $r$-vector norms provided that the local ball is not covered by historical accounts. This serves as the inductive steps in excluding the possibility (ii) in Proposition \ref{P:BNY21}.

    We first fix some notations. Let $\nu$ be a nonzero, compactly supported, $\alpha$-Ahlfors regular, absolutely decaying measure on $\RR^{m\times n}$, and let $\alpha_0>0$ be given as in Proposition \ref{P:BNY21}. Let $u\in\NN$ with $u>\mu^{-1}$, and let $\kappa,\eta>0$ satisfy \begin{equation}\label{E:defkappa}
        0<\kappa<\min\left\{
        \frac{\mu u-1}{2u(\lambda\min\{m,n\}-\mu)},
        \frac{\lambda u-1}{4}
        \right\},
    \end{equation}
    and \begin{equation}\label{E:defeta}
        0<\eta<\min\left\{\frac{\kappa}{4d(d-1)(\kappa+1)},
        \frac{\lambda u\alpha}{\alpha_0}
        \right\}.
    \end{equation}
    Put \begin{equation}\label{E:defalphap}
        \alpha':=\alpha-\frac{\eta\alpha_0}{2\lambda u}.
    \end{equation}

    Let $K=\supp(\nu)$ and consider the $\alpha'$-Cantor potential game on $K$. Let Bob start with a parameter $\beta\in (0,1)$ and a game ball $K\cap B_0$. Here $B_0$ is a ball in $\RR^{m\times n}$ with center in $K$ and radius $\rho_0\in (0,\min\left\{\frac{r_0}{6},\frac{1}{2}\right\})$, where $r_0>0$ is smaller than those in Definition \ref{D:Ahlfors} and Proposition \ref{P:BNY21}. Write $$
    b:=\beta^{-\frac{1}{\lambda}}>1\quad \text{and}\quad \rho_q:=\rho_0\beta^q\quad(q\in\NN_0). 
    $$ 
    Let $c_0>0$ be as in Proposition \ref{P:global-control}. We shall need another constant $c_1>0$ as follows.
    
    \begin{lemma}\label{L:Minkowski}
        There exists a constant $c_1>0$ with the following property. Let $1\leq r<d$ and $\Lambda$ a lattice of $\RR^d$ with covolume $1$. If there exists a rank-$r$ subgroup of $\Lambda$ with covolume $\Delta$ in its real linear span, then $\lambda_1(\Lambda)\leq c_1\Delta^{1/r}$.
    \end{lemma}
    \begin{proof}
        This is a direct consequence of Minkowski's convex body theorem \cite{Ca}.
    \end{proof}

    Now we choose a sufficiently large $N_0\in\NN$ with $N_0\geq 2\kappa^{-1}$ such that for any $N\geq N_0$, \begin{subequations}\label{E:defN_0}
        \begin{align}
            c_0(2\rho_0)^{\min\{m,n\}}\cdot b^{\frac{1}{2}(\mu u-1)N}&\geq 1;\label{E:defN-1}\\
            b^{2\kappa N}\rho_0&\geq 1;\label{E:defN-1.5}\\
            b^{\eta(1+\kappa/2)N}&\geq \sqrt{2}c_1;\label{E:defN-2}\\
            2b^{\eta N}&\geq 3.\label{E:defN-3}
        \end{align}
    \end{subequations}
    The following result is the key proposition in establishing the local estimate.

    \begin{proposition}\label{P:local-control}
        Let $q\in\NN_0$, $B_q$ any ball in $\RR^{m\times n}$ with radius $\rho_q$, and $k,\ell\in\NN_0$ with \begin{equation}\label{E:k+l}
        k+\ell\geq N_0.\end{equation}
        Suppose that there exists some $1\leq r<d$ and a nonzero decomposable $r$-vector $\xi\in\bigwedge^r\ZZ^d$ such that for $t:=q+u(k+\ell)$, \begin{equation}\label{E:local-1}
            \sup_{X\in 2B_q}\left\|\bigwedge^r(D_t^{(k,\ell)}u_X)\xi\right\|<1.
        \end{equation}
        Write $h:=\lfloor \kappa(k+\ell)\rfloor>0$. Then one has $q>uh$; moreover, for any $X\in B_q$,  
        \begin{equation}\label{E:local-control}
            \min\left\{\lambda_1\left(D^{(k+h,\ell)}_tu_X\ZZ^d\right),\lambda_1\left(D^{(k,\ell+h)}_tu_X\ZZ^d\right)\right\}<b^{-\eta(k+\ell+h)}.
        \end{equation}    
        \end{proposition}

    \subsection{Alternative for shrinking norms} In this subsection, we prove some basic facts in exterior algebra, and make a useful observation concerning an alternative for shrinking vector norms under two different diagonal actions. This explains the mechanism behind Proposition \ref{P:local-control}. 

    For a vector $\bff$ in $\RR^d$, we define the contraction by $\bff$ as a linear operator $$
    \iota_{\bff}:\bigwedge^r\RR^d\to\bigwedge^{r-1}\RR^d,\qquad \bv_1\wedge\cdots\wedge\bv_r\mapsto \sum_{j=1}^r(-1)^{j-1}\la \bff,\bv_j\ra\cdot \bv_1\wedge\cdots\wedge\widehat{\bv_j}\wedge\cdots\wedge\bv_r.
    $$
    It is straightforward to check that: \begin{lemma}\label{L:ext}
        Let $\be,\bff$ be two orthonormal vectors in $\RR^d$. For a nonzero decomposable $r$-vector $\xi=\bv_1\wedge\cdots\wedge\bv_r\in \bigwedge^r\RR^d$, write $W=\Span_{\RR}(\bv_1,\cdots,\bv_r)$ and $P_W:\RR^d\to W$ to be the orthogonal projection onto $W$. Then one has 
            \begin{align}              \bigwedge^r(I_d+x\be\cdot\bff^{\top})\xi&=\xi+x(\be\wedge\iota_{\bff}\xi),\quad \forall x\in\RR;\label{E:ext-1}\\
            \|\be\wedge\iota_{\bff}\xi\|^2&=\|\xi\|^2(|\la P_W\be,\bff\ra|^2+\|P_{W^{\perp}}\be\|^2\|P_W\bff\|^2).\label{E:ext-2}
            \end{align}
    \end{lemma}
    \begin{proof}
        The first identity \eqref{E:ext-1} follows directly from the definition. For the second identity \eqref{E:ext-2}, we may assume that $\{\bv_1,\cdots,\bv_r\}$ forms an orthogonal basis of $W$ by the Gram-Schmidt process. We decompose the $r$-vector $\be\wedge\iota_{\bff}\xi$ into an orthogonal sum: $$
        \be\wedge\iota_{\bff}\xi=P_W\be\wedge\iota_{\bff}\xi+P_{W^{\perp}}\be\wedge\iota_{\bff}\xi.
        $$
        The first term equals $\la P_W\be,\bff\ra\cdot\xi$, and the second term is a wedge of two orthogonal vectors. So the norm identity follows.
    \end{proof}

    For two orthonormal vectors $\be,\bff$ in $\RR^d$, we introduce one-parameter diagonal subgroups $\{D_{\bff}^+(t):t\in\RR\}$ and $\{D_{\be}^-(t):t\in\RR\}$ in $\SL_d(\RR)$ by their natural actions on $\RR^d$:\begin{subequations}
    \begin{align}
        D_{\bff}^+(t)\bff=e^{(d-1)t/d}\bff,&\qquad D_{\bff}^+(t)\big|_{\bff^{\perp}}=e^{-t/d}\id_{\bff^{\perp}};\label{E:Df+}\\
        D_{\be}^-(t)\be=e^{-(d-1)t/d}\be,&\qquad D_{\be}^-(t)\big|_{\be^{\perp}}=e^{t/d}\id_{\be^{\perp}}.\label{E:De-}
    \end{align}
    \end{subequations}
    The following lemma shows the alternative for shrinking norms of the $r$-vector $\zeta$ under the two diagonal actions, provided that its spanning subspace $W$ is close enough to both directions $\be$ and $\bff^{\top}$.
    
    \begin{lemma}\label{L:dicho}
Let $s\geq 0$, $T\geq e^{2s}$, $1\leq r<d$, and $\zeta$ a nonzero decomposable $r$-vector with \begin{equation}\label{E:dicho-1}
\sup_{|x|\leq T}\left\|\bigwedge^r(I_d+x\be\cdot\bff^{\top})\zeta\right\|\leq 1.
\end{equation}
Then one has \begin{equation}\label{E:dicho}
\min\left\{
e^{rs/d}\cdot\left\|\bigwedge^rD_{\bff}^+(s)(\zeta)\right\|,\;
e^{(d-r)s/d}\cdot\left\|\bigwedge^rD_{\be}^-(s)(\zeta)\right\|
\right\}
\leq \sqrt{2}.
\end{equation}
\end{lemma}
\begin{proof}
In view of \eqref{E:ext-1}, letting $x=0$ in \eqref{E:dicho-1} gives $\|\zeta\|\leq 1$; letting $x=\pm T$ in \eqref{E:dicho-1} gives $$        \|\zeta+T(\be\wedge\iota_{\bff}\zeta)\|\leq 1\quad \text{and}\quad \|\zeta-T(\be\wedge\iota_{\bff}\zeta)\|\leq 1.
       $$
Then the triangle inequality gives $\|\be\wedge\iota_{\bff}\zeta\|\leq T^{-1}$. By \eqref{E:ext-2}, one has $$        \|\zeta\|\|P_W\bff\|\|P_{W^{\perp}}\be\|\leq \|\be\wedge\iota_{\bff}\zeta\|\leq T^{-1}.
       $$
In particular, it follows that $$        \min\left\{\|P_W\bff\|,\;\|P_{W^{\perp}}\be\|
        \right\}\leq (T\|\zeta\|)^{-1/2}.
       $$
By separating the standard exterior coordinates of $\zeta$ according to whether their basis wedge contains $\bff$, we see from \eqref{E:Df+} that \begin{equation*}
\left\|\bigwedge^rD_{\bff}^+(s)(\zeta)\right\|^2=e^{-2rs/d}\|\zeta\|^2(1+(e^{2s}-1)\|P_W\bff\|^2).
\end{equation*}
Similarly, we see from \eqref{E:De-} that \begin{equation*}
\left\|\bigwedge^rD_{\be}^-(s)(\zeta)\right\|^2=e^{-2(d-r)s/d}\|\zeta\|^2(1+(e^{2s}-1)\|P_{W^{\perp}}\be\|^2).
\end{equation*}
Therefore, the left-hand side of \eqref{E:dicho} equals $$        \|\zeta\|\cdot (1+(e^{2s}-1)\min\{\|P_W\bff\|,\;\|P_{W^{\perp}}\be\|\}^2)^{1/2}
        \leq \|\zeta\|(1+(T-1)(T\|\zeta\|)^{-1})^{1/2}\leq \sqrt{2}.
       $$
This is as desired.
\end{proof}

    \subsection{Proof of Proposition \ref{P:local-control}} This subsection is devoted to proving Proposition \ref{P:local-control}.
    \begin{proof}[Proof of Proposition \ref{P:local-control}]
     We first show that $q>uh$. By Proposition \ref{P:global-control}, one has \begin{equation}\label{E:global-1}
         \begin{aligned}
             \sup_{X\in 2B_q}\left\|\bigwedge^r(D_t^{(k,\ell)}u_X)\xi\right\|
     &\geq c_0b^{\mu t-(k+\ell)}\cdot (2\rho_q)^{\min\{m,n\}}\\
     &=c_0(2\rho_0)^{\min\{m,n\}}\cdot b^{(\mu-\lambda\min\{m,n\})q+(\mu u-1)(k+\ell)}.
         \end{aligned}
     \end{equation}
     Suppose to the contrary that $q\leq uh$. Since $h\leq \kappa(k+\ell)$ and $\mu<\lambda\min\{m,n\}$, we see that the exponent of $b$ is at least $$
     ((\mu-\lambda\min\{m,n\})u\kappa+(\mu u-1))(k+\ell)\overset{\eqref{E:defkappa}}{>}\frac{1}{2}(\mu u-1)(k+\ell).
     $$
     So \eqref{E:defN-1} and \eqref{E:k+l} together imply that the right-hand side of \eqref{E:global-1} is at least $1$. This contradicts \eqref{E:local-1}.

     We next fix any $X_0\in B_q$. Let $E_{i_*,j_*}$ denote the $m\times n$ matrix whose $(i_*,j_*)$-entry is $1$ and other entries are zero. Then for any $|\rho|\leq \rho_q$, one has $X_0+\rho E_{i_*,j_*}\in 2B_q$. Note that $$
     D_t^{(k,\ell)}u_{X_0+\rho E_{i_*,j_*}}=
     (I_d+x\be_{i_*}\cdot \bff_{j_*}^{\top})\cdot D_t^{(k,\ell)}u_{X_0},\quad\text{where}\quad x=b^{\lambda t-(k+\ell)}\rho. 
     $$
     Write
     $$
     T:=b^{\lambda t-(k+\ell)}\rho_q,\quad\text{and}\quad 
     \zeta:=\bigwedge^r(D_t^{(k,\ell)}u_{X_0})\xi.
     $$
     Then it follows from \eqref{E:local-1} that $$
     1>\sup_{|\rho|\leq \rho_q}\left\|
     \bigwedge^r(D_t^{(k,\ell)}u_{X_0+\rho E_{i_*,j_*}})\xi
     \right\|
     =\sup_{|x|\leq T}\left\|\bigwedge^r(I_d+x\be_{i_*}\cdot \bff_{j_*}^{\top})\zeta\right\|.
     $$
     Note that (using \eqref{E:defkappa} and \eqref{E:defN-1.5}) $$
     T=b^{(\lambda u-1)(k+\ell)}\rho_0\geq b^{2h}=e^{2h\log{b}}.$$
     Then we see from Lemma \ref{L:dicho} that $$
     \min\left\{
            b^{rh/d}\cdot\left\|\bigwedge^rD_{\bff_{j_*}}^+(h\log{b})(\zeta)\right\|,\;
            b^{(d-r)h/d}\cdot\left\|\bigwedge^rD_{\be_{i_*}}^-(h\log{b})(\zeta)\right\|
            \right\}
            \leq \sqrt{2}.
     $$
     Note that $$
     D_{\bff_{j_*}}^+(h\log{b})\cdot D_t^{(k,\ell)}=D_t^{(k,\ell+h)}\quad \text{and}\quad D_{\be_{i_*}}^-(h\log{b})\cdot D_t^{(k,\ell)}=D_t^{(k+h,\ell)}.
     $$
     Therefore, we conclude that $$
     \min\left\{
            b^{rh/d}\cdot\left\|\bigwedge^r(D_t^{(k,\ell+h)}u_{X_0})\xi\right\|,\;
            b^{(d-r)h/d}\cdot\left\|\bigwedge^r(D_t^{(k+h,\ell)}u_{X_0})\xi\right\|
            \right\}
            \leq \sqrt{2}.
     $$

     Finally, by Lemma \ref{L:Minkowski}, we have for any diagonal matrix $D$, $$
     \lambda_1(Du_{X_0}\ZZ^d)\leq c_1\cdot\left\|\bigwedge^r(Du_{X_0})\xi\right\|^{1/r}.
     $$
     It follows from \eqref{E:defeta} that $$2\eta(k+\ell+h)<\min\left\{\frac{h}{d},\;\frac{(d-r)h}{rd}\right\}.$$
    Combining these estimates together with \eqref{E:defN-2} gives \eqref{E:local-control}.
     \end{proof}

     \subsection{Measure control of dangerous sets}\label{SS:measure-control}
     In this subsection, we deduce the local estimate for the $r$-vector norms from Proposition \ref{P:local-control}. Combined with Proposition \ref{P:BNY21}, this gives a local control over the measure of those dangerous sets.

     For $q,k,\ell\in\NN_0$ with $k+\ell\geq N_0$, we define the \textsl{raw $(q;k,\ell)$-dangerous set} as \begin{equation}
         \cD_q^{(k,\ell)}:=\left\{
         X\in\RR^{m\times n}:\lambda_1\left(D^{(k,\ell)}_{q+u(k+\ell)}u_X\ZZ^d\right)<b^{-\eta(k+\ell)}
         \right\}.
     \end{equation}
     Given Bob's $q$-th game ball $K\cap B_q$, we define the \textsl{current $(q;k,\ell)$-dangerous set} as \begin{equation}
        \widetilde{\cD}_q^{(k,\ell)}:=(K\cap B_q\cap \cD_q^{(k,\ell)})\setminus \bigcup_{q'<q}\bigcup_{k'+\ell'\geq N_0}\cD_{q'}^{(k',\ell')}.
     \end{equation}
     We focus on the nontrivial case when some current dangerous set is nonempty.
     
     \begin{corollary}\label{C:local-control}
         Let $C_0,\alpha_0>0$ be given as in Proposition \ref{P:BNY21}, and $q\in\NN_0$. Suppose that \begin{equation}\label{E:history}
             (K\cap B_q)\setminus \bigcup_{0\leq q'<q}\bigcup_{k'+\ell'\geq N_0}\cD_{q'}^{(k',\ell')}\neq \varnothing.
         \end{equation}
         Then for any $(k,\ell)\in \NN_0^2$ with $k+\ell\geq N_0$, any $1\leq r<d$, and any nonzero decomposable $r$-vector $\xi\in \bigwedge^r\ZZ^d$, one has \begin{equation}\label{E:local-control-C}
             \sup_{X\in 2B_q}\left\|
             \bigwedge^r(D^{(k,\ell)}_{q+u(k+\ell)}u_X)\xi
             \right\|\geq 1.
         \end{equation} 
         In particular, for any $0<\epsilon\leq 1$, one has \begin{equation}\label{E:measure-control}
             \nu\left\{
             X\in 2B_q:\lambda_1(D^{(k,\ell)}_{q+u(k+\ell)}u_X\ZZ^d)<\epsilon
             \right\}\leq C_0\epsilon^{\alpha_0}\cdot\nu(6B_q).
         \end{equation}
     \end{corollary}
     \begin{proof}
     Suppose that \eqref{E:local-control-C} fails for some $k,\ell,r,\xi$. Then for $h=\lfloor \kappa(k+\ell)\rfloor>0$, Proposition \ref{P:local-control} gives $q':=q-uh\in (0,q)$, and for any $X\in B_q$, \eqref{E:local-control} holds. This implies that $$
     B_q\subseteq \cD^{(k+h,\ell)}_{q'}\cup \cD^{(k,\ell+h)}_{q'},
     $$
     which contradicts \eqref{E:history}. Applying Proposition \ref{P:BNY21} with $B=2B_q$, $D=D^{(k,\ell)}_{q+u(k+\ell)}$ and $\rho=1$ gives the final statement.
     \end{proof}

     \section{Proofs of Theorems \ref{T:main} and \ref{T:dyn}}\label{S:Proofs}
     In this section, we prove our first main result Theorem \ref{T:main} and extend it to weights with zero components. The dynamical statement Theorem \ref{T:dyn} is actually equivalent to this extended result. 

     \subsection{Cantor potential and winning strategy} In this subsection, we verify that the measure control of dangerous sets in Section \ref{SS:measure-control} provides a legal collection of balls for Alice to delete.

     \begin{proposition}\label{P:Cantor}
         Let $A_0,\alpha$ be as in Definition \ref{D:Ahlfors}, and $C_0,\alpha_0>0$ as in Proposition \ref{P:BNY21}. Set \begin{equation}
             M_0=A_0^2C_0(3/2)^{\alpha_0}12^{\alpha}.
         \end{equation}
    For $q,k,\ell\in\NN_0$ with $k+\ell\geq N_0$, the set $\widetilde{\cD}_q^{(k,\ell)}$ admits a cover $\cA^{(k,\ell)}_q$ consisting of balls with centers in $K$ and radius $\rho_{q+u(k+\ell)}$, whose cardinality is at most $M_0b^{(\lambda u\alpha-\eta\alpha_0)(k+\ell)}$.
    \end{proposition}
    \begin{proof}
    If $\widetilde{\cD}_q^{(k,\ell)}=\varnothing$, then there is nothing to cover. Otherwise, write $$
    D=D_{q+u(k+\ell)}^{(k,\ell)},\quad \widetilde{\cD}=\widetilde{\cD}_q^{(k,\ell)},\quad \epsilon=b^{-\eta(k+\ell)},\quad r=\rho_{q+u(k+\ell)}/2.
    $$
    We claim that: for the ambient $r$-neighborhood of $\widetilde{\cD}$ in $\RR^{m\times n}$, one has \begin{equation}\label{E:danger-contain}
        \widetilde{\cD}^{(r)}\subseteq \left\{
        Y\in 2B_q:\lambda_1(Du_Y\ZZ^d)<3\epsilon/2
        \right\}.
    \end{equation}

    Indeed, let $Y\in \widetilde{\cD}^{(r)}$. Then there is some $X\in \widetilde{\cD}$ with $\|Y-X\|\leq r$. We see that $$
    Du_Y=Du_{Y-X}D^{-1}\cdot Du_X,
    $$
    where $$
    \|Du_{Y-X}D^{-1}\|\leq 1+b^{\lambda(q+u(k+\ell))}\|Y-X\|\leq 1+\rho_0/2<3/2.
    $$
    Since $\lambda_1(Du_X\ZZ^d)<\epsilon$, we conclude that $\lambda_1(Du_Y\ZZ^d)<3\epsilon/2$. This verifies \eqref{E:danger-contain}.

    By Corollary \ref{C:local-control} (using \eqref{E:defN-3}), one has $$
    \nu(\widetilde{\cD}^{(r)})\overset{\eqref{E:danger-contain}}{\leq} \nu\left\{
        Y\in 2B_q:\lambda_1(Du_Y\ZZ^d)<3\epsilon/2
        \right\}
        \overset{\eqref{E:measure-control}}{\leq }C_0(3\epsilon/2)^{\alpha_0}\cdot\nu(6B_q).
    $$
    By Ahlfors regularity of $\nu$, one has $\nu(6B_q)\leq A_0(6\rho_q)^{\alpha}$. Therefore, Lemma \ref{L:BNY-cover} gives a cover of $\widetilde{\cD}$ consisting of balls with centers in $K$ and radius $2r=\rho_{q+u(k+\ell)}$, whose cardinality is at most $$
    A_0\nu(\widetilde{\cD}^{(r)})\cdot r^{-\alpha}\leq A_0^2C_0(3\epsilon/2)^{\alpha_0}(12\beta^{-u(k+\ell)})^{\alpha}=M_0b^{(\lambda u\alpha-\eta\alpha_0)(k+\ell)}.
    $$
    This is as desired.
    \end{proof}

    Now we enlarge $N_0\in\NN$ in Section \ref{S:local-estimate} so that for any $N\geq N_0$, \begin{equation}\label{E:defN-4}
        b^{\eta\alpha_0N/2}\geq (N+1)M_0.
    \end{equation}
    We are ready to define Alice's strategy of deleting collections of balls. For $q\in\NN_0$, let Alice choose her collections $\cA_{q+1,i}\,(i\geq 0)$ as \begin{equation}
        \cA_{q+1,i}:=\begin{cases}
           \bigcup_{(k,\ell):k+\ell=N}\cA^{(k,\ell)}_q, &\text{if }i=uN-1\text{ for some }N\geq N_0;\\
            \varnothing, &\text{otherwise}
        \end{cases}
    \end{equation}
    By Proposition \ref{P:Cantor}, the balls in $\cA_{q+1,i}\neq\varnothing$ have radius $\rho_{q+uN}=\rho_{q+1+i}$, and by \eqref{E:defN-4}\eqref{E:defalphap}, $$
    |\cA_{q+1,i}|\leq \sum_{(k,\ell):k+\ell=N}|\cA^{(k,\ell)}_q|\leq (N+1)M_0b^{(\lambda u\alpha-\eta\alpha_0)N}\leq \beta^{-\alpha'(i+1)}.
    $$
    So Alice's choices are valid.

    Finally, we show that Alice's strategy above is a winning strategy with the target set $\Bad(\br;\bs)$, and hence complete the proof of Theorem \ref{T:main}.

    \begin{proof}[Proof of Theorem \ref{T:main}] Let $\nu$ be a nonzero, compactly supported, $\alpha$-Ahlfors regular, absolutely decaying measure on $\RR^{m\times n}$. Write $K=\supp(\nu)$. We choose $\alpha'\in (0,\alpha)$ as in \eqref{E:defalphap} and let the $\alpha'$-Cantor potential game on $K$ go as in Section \ref{S:local-estimate}.
    
    We first show that: if Alice doesn't win by default, then the outcome of the game avoids all raw dangerous sets $\cD_{q}^{(k,\ell)}$ for $q,k,\ell\in\NN_0$ with $k+\ell\geq N_0$. In fact, let $X_*\in K\cap\bigcap_{q=0}^{+\infty}B_q$. Suppose that $X_*\in \cD^{(k,\ell)}_q$ for some $q,k,\ell$. Then we may choose $q$ minimal so that $$
    X_*\in \cD_q^{(k,\ell)}\setminus \bigcup_{q'<q}\bigcup_{k'+\ell'\geq N_0}\cD_{q'}^{(k',\ell')}.
    $$ 
    This implies that $X_*\in \widetilde{\cD}^{(k,\ell)}_q$. So $X_*$ is contained in some ball in the collection $\cA^{(k,\ell)}_q$, which doesn't meet $K\cap B_{q+u(k+\ell)}$. This is a contradiction.

    We next show that: if some point $X_*$ avoids all raw dangerous sets $\cD_{q}^{(k,\ell)}$ for $q,k,\ell\in\NN_0$ with $k+\ell\geq N_0$, then $X_*\in \Bad(\br;\bs)$. In fact, for all $q,k,\ell$, one has $$   \lambda_1\left(D^{(k,\ell)}_{q+u(k+\ell)}u_{X_*}\ZZ^d\right)\geq b^{-\eta(k+\ell)}.
    $$
    In particular, we take $(k,\ell)=(N_0,0)$, and use $$D^{(N_0,0)}_{0}\cdot g^{(\br;\,\bs)}_{q+uN_0}=D^{(N_0,0)}_{q+uN_0}$$
    to get for all $q\in\NN_0$, \begin{equation*}
       \lambda_1\left(g^{(\br;\,\bs)}_{q+uN_0}u_{X_*}\ZZ^d\right)
    \geq \|D^{(N_0,0)}_{0}\|^{-1}\cdot\lambda_1\left(D^{(N_0,0)}_{q+uN_0}u_{X_*}\ZZ^d\right)
    \geq b^{-(\eta+1/d)N_0}.
    \end{equation*}
    Moreover, for $0\leq q<uN_0$,  one has $$
    \lambda_1(g^{(\br;\,\bs)}_qu_{X_*}\ZZ^d)\geq \|g^{(\br;\,\bs)}_{uN_0}\|^{-1}\cdot \lambda_1\left(g^{(\br;\,\bs)}_{q+uN_0}u_{X_*}\ZZ^d\right)\geq b^{-(r_{i_*}u+\eta+1/d)N_0}.
    $$
    So we conclude that $$
    \inf_{q\in\NN_0}\lambda_1(g^{(\br;\,\bs)}_qu_{X_*}\ZZ^d)>0.
    $$
    Then it follows from Dani-Kleinbock correspondence that $X_*\in\Bad(\br;\bs)$. 

    Therefore, we have shown that the target set $K\cap\Bad(\br;\bs)$ is $\alpha'$-Cantor winning. Then Theorem \ref{T:BNY21} (ii) shows that $K\cap\Bad(\br;\bs)\neq\varnothing$. By Theorem \ref{T:BNY21} (i), the set $\Bad(\br;\bs)$ is HAW.
    \end{proof}

    \subsection{Extension to weights with zero components}\label{SS:zeroweight} In this subsection, we first extend the definition of $\Bad(\br;\bs)$ to weights with zero components, and then prove that the resulting set is also HAW.
    
    To avoid any possible ambiguity, we don't adopt the imprecise notations like $1/0=\infty$. Instead, we redefine $\Bad(\br;\bs)$ for all pairs of weights possibly with zero components, and show that this coincides with the definition \eqref{E:Bad} using quasi-norms.

    For $a\in\NN$, let $\overline{\cW_a}$ denote the set of all $a$-dimensional weights possibly with zero components: \begin{equation}
        \overline{\cW_a}:=\left\{
        \bw=(w_1,\cdots,w_a)\in\RR^a:w_1,\cdots,w_a\geq 0,\;\sum_{i=1}^aw_i=1
        \right\}.
    \end{equation}
    It is the topological closure of $\cW_a$ in $\RR^a$. For $\bw\in \overline{\cW_a}$, write \begin{equation}
        J_{\bw}:=\{1\leq i\leq a:w_i\neq 0\},\quad \text{and}\quad \RR^{J_{\bw}}:=\{\bv\in\RR^a:v_i=0,\;\forall i\notin J_{\bw}\}.
    \end{equation}

    Let $m,n\in\NN$. For $(\br;\bs)\in\overline{\cW_m}\times \overline{\cW_n}$, we say an $m\times n$ real matrix $X=(X_{ij})$ is $\overline{(\br;\bs)}$-\textsl{badly approximable}, if \begin{equation}\label{E:Bad0}
        \inf_{\substack{\bp\in\ZZ^m\\\bq\in \ZZ^{J_{\bs}}\setminus\{0\}}}\max_{1\leq i\leq m}
        \left(\max_{j\in J_{\bs}}|q_j|^{1/s_j}\right)^{r_i}\cdot \left|\sum_{j\in J_{\bs}}X_{ij}q_j+p_i\right|>0.
    \end{equation}
    The set of $\overline{(\br;\bs)}$-badly approximable matrices is denoted by $\Bad\overline{(\br;\bs)}$. It is straightforward to see that: when $(\br;\bs)\in\cW_m\times \cW_n$, one has \begin{equation}
        \Bad\overline{(\br;\bs)}=\Bad(\br;\bs).
    \end{equation}
    The main result of this subsection is the following:

    \begin{proposition}\label{P:extend}
        For any $(\br;\bs)\in \overline{\cW_m}\times \overline{\cW_n}$, the set $\Bad\overline{(\br;\bs)}$ is HAW.
    \end{proposition}
    \begin{proof}
        Put $I:=J_{\br}$ and $J:=J_{\bs}$, and write $m':=|I|$ and $n':=|J|$. Since the components of each weight sum to $1$, both $I$ and $J$ are nonempty. The restricted weights
        $$\br':=(r_i)_{i\in I}\in\cW_{m'},\qquad
        \bs':=(s_j)_{j\in J}\in\cW_{n'}$$
        have positive components. Let
        $$\pi:\RR^{m\times n}\longrightarrow\RR^{m'\times n'},\qquad
        \pi(X)=(X_{ij})_{i\in I,\,j\in J}$$
        be the coordinate projection. Then it is straightforward to see that
        $$\pi^{-1}\bigl(\Bad(\br';\bs')\bigr)
        \subseteq\Bad\overline{(\br;\bs)}.$$
        The conclusion follows from Theorem \ref{T:main} and \cite[Lemma 2.2 (v)]{AGK}.   
    \end{proof}

    \subsection{Proof of Theorem \ref{T:dyn}} In this subsection, we first verify the dynamical interpretation of the set $\Bad\overline{(\br;\bs)}$ in the language of bounded orbits, and then deduce Theorem \ref{T:dyn} from Proposition \ref{P:extend}. We start with a simple fact from Dani-Kleinbock correspondence. 

    \begin{lemma}\label{L:Dani-Kleinbock}
        An $m\times n$ real matrix $X$ is $\overline{(\br;\bs)}$-badly approximable if and only if the trajectory $\{g_t^{(\br;\,\bs)}u_X\Gamma:t\geq 0\}$ is bounded in $G/\Gamma=\SL_d(\RR)/\SL_d(\ZZ)$.
    \end{lemma}
    \begin{proof}
        Write $J:=J_{\bs}$ and
        $$Q(\bq):=\max_{j\in J}|q_j|^{1/s_j}
        \quad\text{for }\bq\in\ZZ^J\setminus\{0\}.$$
        In particular, $Q(\bq)\geq 1$. We use the maximum norm $\|\cdot\|_{\infty}$ on $\RR^d$; its equivalence with the Euclidean norm shows that this does not change the boundedness criterion. By Mahler's compactness criterion, the trajectory in the statement is bounded if and only if there exists $\delta>0$ such that
        \begin{equation}\label{E:Mahler}
            \left\|g_t^{(\br;\,\bs)}u_X(\bp,\bq)^{\top}\right\|_{\infty}\geq\delta
        \quad\text{for all }t\geq 0\text{ and }(\bp,\bq)^{\top}\in\ZZ^d\setminus\{0\}.
        \end{equation}

        First suppose that $X\in\Bad\overline{(\br;\bs)}$. By \eqref{E:Bad0}, there is $c>0$ such that
        $$\max_{1\leq i\leq m}Q(\bq)^{r_i}|(X\bq+\bp)_i|\geq c
        \quad\text{for all }\bp\in\ZZ^m,\ \bq\in\ZZ^J\setminus\{0\}.$$
        Fix $t\geq0$ and $(\bp,\bq)^{\top}\in\ZZ^d\setminus\{0\}$. If $q_j\neq0$ for some $j\notin J$, then the corresponding coordinate of $g_t^{(\br;\,\bs)}u_X(\bp,\bq)^{\top}$ equals $q_j$ and has absolute value at least $1$. If $\bq=0$, then some $p_i\neq0$, and the corresponding coordinate has absolute value $b^{r_it}|p_i|\geq1$.

        It remains to consider $\bq\in\ZZ^J\setminus\{0\}$. If $b^t<Q(\bq)$, choose $j\in J$ with $|q_j|=Q(\bq)^{s_j}$. Then
        $$b^{-s_jt}|q_j|=(Q(\bq)/b^t)^{s_j}>1.$$
        If $b^t\geq Q(\bq)$, one has
        $$\max_{1\leq i\leq m}b^{r_it}|(X\bq+\bp)_i|
        \geq\max_{1\leq i\leq m}Q(\bq)^{r_i}|(X\bq+\bp)_i|\geq c.$$
        Thus the maximum norm in \eqref{E:Mahler} is always at least $\min\{1,c\}>0$.

        Conversely, suppose that there exists some $\delta\in(0,1)$ satisfying \eqref{E:Mahler}. Set
        $$M:=(2/\delta)^{1/\min_{j\in J}s_j}>1.$$
        Fix $\bp\in\ZZ^m$ and $\bq\in\ZZ^J\setminus\{0\}$, and choose $t\geq0$ such that $b^t=MQ(\bq)$. For every $j\in J$,
        $$b^{-s_jt}|q_j|\leq M^{-s_j}\leq\delta/2,$$
        while all coordinates with $j\notin J$ vanish. So \eqref{E:Mahler} implies that, for some $1\leq i\leq m$,
        $$b^{r_it}|(X\bq+\bp)_i|\geq\delta.$$
        Consequently,
        $$\max_{1\leq i\leq m}Q(\bq)^{r_i}|(X\bq+\bp)_i|
        \geq\delta M^{-r_i}\geq\delta M^{-\max_{1\leq i\leq m}r_i}>0.$$
        This lower bound is independent of $\bp$ and $\bq$, so \eqref{E:Bad0} holds. The proof is complete.
    \end{proof}

    \begin{proof}[Proof of Theorem \ref{T:dyn}] Let $F^+$ be a nontrivial ray in $A^+$, and write $F^+=(F^{(\br;\bs)})^+$ for some $(\br;\bs)\in \overline{\cW_m}\times \overline{\cW_n}$. Since the map $$
    \RR^{m\times n}\to \TT_0,\quad X\mapsto u_X\Gamma
    $$
    is a local diffeomorphism, by Lemma \ref{L:HAW-Ra} (iii), it suffices to show that the set $$
    \{X\in\RR^{m\times n}:(F^{(\br;\bs)})^+u_X\Gamma\text{ is bounded}\}
    $$
    is HAW. In view of Lemma \ref{L:Dani-Kleinbock}, the above set equals $\Bad\overline{(\br;\bs)}$. Then the conclusion follows from Proposition \ref{P:extend}.
    \end{proof}

		\section*{Acknowledgment} Part of this project was initiated during the author's self-funded visit to Brandeis University in 2024--2025. The author is grateful to Dmitry Kleinbock for his generous hospitality and warm encouragement throughout this period. Thanks are also due to Jinpeng An, Lifan Guan, Bohan Yang, and Runlin Zhang for their helpful discussions.
		
		\section*{Declaration of AI Use}
        The current version of this paper is developed with the help of ChatGPT, especially in simplifying arguments from Beresnevich-Nesharim-Yang's papers \cite{BNY21,BNY22}. The author has independently verified all mathematical arguments and takes full responsibility for the manuscript.

	\end{document}